\RequirePackage{ifpdf}
\ifpdf 
\documentclass[pdftex]{sigma}
\else
\documentclass{sigma}
\fi

\def\R{\mathbb{R}}
\def\C{\mathbb{C}}

\def\k{\kappa}

\def\g{\mathfrak{g}}

\def\p{\mathfrak{p}}

\def\Abdle{\mathcal{A}}

\def\Gbdle{\mathcal{G}}

\def\VF{\mathfrak{X}}
\def\Lie{\mathcal{L}}
\def\:{\lrcorner}
\def\#{\sharp}

\def\tens{\otimes}

\def\dsum{\oplus}

\def\prod{\times}

\def\isom{\cong}

\theoremstyle{plain}
\newtheorem{Theorem}{Theorem}[section]
\newtheorem{Lemma}[Theorem]{Lemma}

\newtheorem{Proposition}[Theorem]{Proposition}

\theoremstyle{definition}
\newtheorem{Definition}{Definition}[section]
\newtheorem{Remark}[Definition]{Remark}

\newtheorem*{TheoremA}{Theorem A}
\newtheorem*{TheoremB}{Theorem B}

\begin{document}

\allowdisplaybreaks

\renewcommand{\PaperNumber}{039}

\FirstPageHeading

\ShortArticleName{Essential Parabolic Structures and Their Inf\/initesimal Automorphisms}

\ArticleName{Essential Parabolic Structures\\ and Their Inf\/initesimal Automorphisms}

\Author{Jesse ALT}

\AuthorNameForHeading{J.~Alt}

\Address{School of Mathematics, University of the Witwatersrand,\\ P O Wits 2050, Johannesburg, South Africa}
\Email{\href{mailto:jesse.alt@wits.ac.za}{jesse.alt@wits.ac.za}}
\URLaddress{\url{http://sites.google.com/site/jmaltmath/}}

\ArticleDates{Received November 02, 2010, in f\/inal form April 11, 2011;  Published online April 14, 2011}

\Abstract{Using the theory of Weyl structures, we give a natural generalization of the notion of essential conformal structures and conformal Killing f\/ields to arbitrary parabolic geometries. We show that a parabolic structure is inessential whenever the automorphism group acts properly on the base space. As a corollary of the generalized Ferrand--Obata theorem proved by C.~Frances, this proves a generalization of the ``Lichn\'{e}rowicz conjecture'' for conformal Riemannian, strictly pseudo-convex CR, and quaternionic/octonionic contact manifolds in positive-def\/inite signature. For an inf\/initesimal automorphism with a~singularity, we give a generalization of the dictionary introduced by Frances for conformal Killing f\/ields, which characterizes (local) essentiality via the so-called holonomy associated to a~singularity of an inf\/initesimal automorphism.}


\Classification{53B05; 53C05; 53C17; 53C24}

\section{Introduction}

\subsection{Motivation from conformal geometry}

Let $(M,c)$ be a smooth, $n$-dimensional semi-Riemannian conformal manifold. For any choice of semi-Riemannian metric $g$ from the equivalence class $c$ def\/ining the conformal structure, we have the obvious inclusion of the group of isometric dif\/feomorphisms of $(M,g)$ in the group of conformal dif\/feomorphisms of $(M,c)$, $\mathrm{Isom}(M,g) \subseteq \mathrm{Conf}(M,c)$. At the inf\/initesimal level of vector f\/ields, we have the corresponding inclusion of Killing f\/ields in the conformal vector f\/ields, $\mathrm{KVF}(M,g) \subseteq \mathrm{CVF}(M,c)$, which is obvious from the def\/initions: $\mathrm{KVF}(M,g) := \{ X \in \VF(M) \, \vert \, \mathcal{L}_X g = 0\}$; and $\mathrm{CVF}(M,c) := \{ X \in \VF(M) \, \vert \, \exists \, \lambda \in C^{\infty}(M) \, \mathrm{s.t.} \, \mathcal{L}_X g = \lambda g   \}$.

A conformal dif\/feomorphism $\varphi \in \mathrm{Conf}(M,c)$ is called essential if $\varphi$ is not an isometry of any metric $g \in c$, and the conformal structure $(M,c)$ is essential if $\mathrm{Isom}(M,g)$ is a proper subgroup of $\mathrm{Conf}(M,c)$ for all representatives $g \in c$. Similarly, a conformal vector f\/ield $X \in \mathrm{CVF}(M,c)$ is called essential if there is no representative $g \in c$ for which $X \in \mathrm{KVF}(M,g)$. It is a fact~-- although not necessarily obvious from the preceding def\/initions -- that there are compact and non-compact essential conformal structures in all dimensions $n \geq 2$ and all signatures $(p,q)$, which moreover admit essential conformal vector f\/ields. The standard compact example is given by the conformal ``M\"obius sphere'' $(S^{p,q},c)$ of any signature $(p,q)$ (also called the Einstein universe~-- these are the conformally f\/lat homogeneous models of conformal geometry, which in Riemannian signature are just the standard $n$-spheres equipped with the conformal class of the round metric), while the standard non-compact example is $\R^{p+q}$ equipped with the conformal class of the f\/lat metric of signature $(p,q)$. In fact, as a result of the following well-known theorem, giving a positive answering to the so-called Lichn\'{e}rowicz conjecture, we know that in Riemannian signature these two examples are the only essential structures:

\begin{TheoremA}[Ferrand--Obata for $n \geq 3$] \label{TheoremA}
\it If $(M,c)$ is an essential Riemannian conformal structure of dimension $n \geq 2$, then it is conformally diffeomorphic to the $n$-dimensional sphere with the round metric, or to $n$-dimensional Euclidean space.
\end{TheoremA}

For compact manifolds, this theorem was proven by M.~Obata and J.~Lelong-Ferrand in the late 1960's and early 1970's. A proof for the non-compact case, announced in 1972 by Alekseevski, was later discovered to be incomplete, and a complete proof was f\/irst given in 1994 by Ferrand (cf.~\cite{Ferrand, FrEss} and references therein). Recently, a corresponding result was proven at the inf\/initesimal level by C.~Frances~\cite{Frances3} (note that this theorem does not simply follow from an application of the Ferrand--Obata theorem, because the conformal vector f\/ields are not assumed to be complete):

\begin{TheoremB}[Frances] \label{TheoremB}
\it Let $(M,c)$ be a conformal Riemannian manifold of dimension $n \geq 3$, endowed with a conformal vector f\/ield~$X$ which vanishes at $x_0 \in M$. Then either: $(1)$ There exists a neighborhood~$U$ of~$x_0$ on which~$X$ is complete, generates a relatively compact flow in $\mathrm{Conf}(U,c)$, and is inessential on~$U$, i.e.\ $X \in \mathrm{KVF}(U,g)$ for some $g \in c_{\vert U}$; or $(2)$ There is a~conformally flat neighborhood $U$ of~$x_0$, and~$X$ is essential on each neighborhood of~$x_0$.
\end{TheoremB}

\subsection{Organization of the text and summary of main results}

One direction of research into how these results do (or do not) generalize to other settings is to consider the analogous questions for \emph{pseudo}-Riemannian metrics, where essential conformal structures turn out to be much more prevalent (cf.~\cite{FrEss} for a survey). The aim of the present text is to introduce natural generalizations of the notion of essential structure, and the corresponding notion at the inf\/initesimal level, to the category of parabolic geometries. After this, we establish a generalization of Theorem~A to a class of geometries which have been called ``rank one parabolic geometries'': conformal Riemannian structures; strictly pseudo-convex CR structures of hypersurface type; positive-def\/inite quaternionic contact structures; and octonionic contact structures (cf.~\cite{Biquard}). In fact, once our general def\/initions have been introduced and some basic properties established, we only have to prove the easy part of this generalized Theorem~A, the dif\/f\/icult part having been taken care of in \cite{Frances}; in the CR case, similar results to \cite{Frances} were obtained in \cite{Schoen} and \cite{Webster}; in the quaternionic contact case, see \cite{IvVass}. Finally, we establish some local properties of essential inf\/initesimal automorphisms, which generalize essential conformal vector f\/ields.

We begin in Section~\ref{section2.1} with a review of relevant tools from parabolic geometry, in particular the notions of Weyl structures introduced in \cite{CS03}, which are then used to generalize the notions of essentiality to arbitrary parabolic geometries, cf.\ Def\/inition~\ref{parabolic essential definition}. Next, we turn to some basic properties of essential automorphisms in Section~\ref{section2.2}, establishing equivalent characterizations which will be required in the proofs of the main results.

In Section~\ref{section3} we establish a basic global result in the general parabolic setting: a parabolic structure is essential only if the action of the automorphism group $\mathrm{Aut}(\Gbdle,\omega)$ on $M$ is non-proper (cf.\ Proposition~\ref{proper implies inessential}, which is a generalization of a result proven in the conformal case by Alekseevski in~\cite{Aleks}). With this, we may apply the main theorem of~\cite{Frances} to prove a Lichn\'{e}rowicz theorem for rank one parabolic geometries, conf\/irming the conjecture formulated in Section~\ref{section2.2} of~\cite{FrEss} for these geometries:

\begin{Theorem} \label{Lich for rank one parabolic} Let $(\Gbdle \rightarrow M,\omega)$ be a regular rank one parabolic geometry, with $M$ connected. If this parabolic structure is essential, then $M$ is geometrically isomorphic to either the compact homogeneous model $G/P$ or the noncompact space $G/P \backslash \{eP\}$.
\end{Theorem}

In Section~\ref{section4} we establish some local properties of essential inf\/initesimal automorphisms which generalize some of the results of~\cite{Frances3}. We begin in Section~\ref{section4.1} by recalling a result characterizing inf\/initesimal automorphisms of arbitrary Cartan geometries $(\Gbdle \rightarrow M,\omega)$ via an identity involving the curvature of $\omega$. This identity was established in~\cite{Cap inf aut} for parabolic geometries and carries over without dif\/f\/iculty to general Cartan geometries. Next, we show that the local study of essential inf\/initesimal automorphisms amounts to studying their singularities, since any inf\/initesimal automorphism of a parabolic geometry is inessential in some neighborhood of any point $x$ such that $X(x) \neq 0$ (cf.\ Proposition~\ref{non-singular points}). Then, we apply the identity reviewed in Section~\ref{section4.1} to prove a generalization of results of~\cite{Frances3}, which give a ``dictionary'' relating essentiality of an inf\/initesimal automorphism near a singularity~$x_0$ to properties of its holonomy~$h^t$ at~$x_0$, a one-parameter subgroup of $P$ which is determined up to conjugacy (cf.\ Def\/inition \ref{holonomy definition}; it should be emphasized that this notion of ``holonomy'' of the singularity of an inf\/initesimal automorphism is distinct from the holonomy of, e.g., a Cartan connection which is common in the literature). The main local result can be stated as:

\begin{Theorem} \label{local essential dictionary} Let $(\Gbdle \rightarrow M,\omega)$ be a parabolic geometry of type $(G,P)$ and $\mathbf{X} \in \mathrm{inf}(\Gbdle,\omega)$ an infinitesimal automorphism with singularity at $x_0 \in M$. Then $\mathbf{X}$ is inessential in some neighborhood $U$ of $x_0$ if and only if the holonomy $h^t$ of $\mathbf{X}$ at $x_0$ is up to conjugacy a subgroup of $\mathrm{Ker}(\lambda) \subset G_0$ $($equivalently, if and only if $\omega(\mathbf{X}(u_0)) \in \mathrm{Ker}(\lambda') \subset \g_0$ for some $u_0 \in \Gbdle_{x_0})$.
\end{Theorem}

Already in conformal geometry, this result is of some interest because it can be used to determine whether a conformal vector f\/ield is locally essential from looking at the adjoint tractor it determines. We expect that the generalization to arbitrary parabolic geometries will be useful in trying to generalize Theorem B to the other rank one parabolic geometries.

\section{Essential automorphisms: basic def\/initions and properties}\label{section2}

\subsection{Background on parabolic geometries and their Weyl structures}\label{section2.1}

Let us begin by recalling the def\/initions of parabolic geometries and their Weyl structures (the latter, introduced by A.~\v{C}ap and J.~Slov\'ak in \cite{CS03}, will be central to our notion of essential parabolic structures). Parabolic geometries are certain types of Cartan geometries, which are very general: given a closed subgroup $P$ of a Lie group $G$, a \textit{Cartan geometry of type $(G,P)$} (or \textit{modelled on the homogeneous space $G/P$}) is given by a principal $P$ bundle $\pi: \Gbdle \rightarrow M$, equipped with a \textit{Cartan connection} $\omega$. That is, $\omega \in \Omega^1(\Gbdle,\g)$ satisf\/ies:
\begin{gather}
R_p^*(\omega) = \mathrm{Ad}\big(p^{-1}\big) \circ \omega, \,\, \mathrm{for} \,\, \mathrm{all} \,\, p \in P; \label{ad-inv}\\
\omega(\tilde{X}) = X, \,\, \mathrm{for} \,\, \mathrm{any} \,\, X \in \p, \, \tilde{X} \,\, \mathrm{its} \,\, \mathrm{fundamental} \,\, \mathrm{vector} \,\, \mathrm{f\/ield} \,\, \mathrm{on} \,\, \Gbdle; \label{fund-vfs}\\
\omega(u) : T_u\Gbdle \rightarrow \g \,\, \mathrm{is} \,\, \mathrm{a} \,\, \mathrm{linear} \,\, \mathrm{isomorphism} \,\, \mathrm{for} \,\, \mathrm{all} \,\,u \in \Gbdle. \label{parallelism}
\end{gather}

A Cartan geometry of type $(G,P)$ is a \emph{parabolic geometry} if $G$ is a real or complex semi-simple Lie group, and $P \subset G$ is a parabolic subgroup as in representation theory~-- at the level of Lie algebras, this means, for $\g$ complex semi-simple, that $\p$ must contain a Borel (maximal solvable) subalgebra of $\g$; for $\g$ real semi-simple, the complexif\/ication $\p(\C)$ must contain a Borel subalgebra of $\g(\C)$. (For a more detailed discussion of the basic properties of parabolic subgroups and parabolic geometries, the reader is referred to \cite{CSbook}. Here we only attempt to cite some of the key facts which are germane to the subsequent text.) In particular, in the parabolic setting the Lie algebra $\g$ of $G$ has an induced $\vert k \vert$-grading for some natural number $k$, so $\g = \g_{-k} \dsum \cdots \dsum \g_k$ with $[\g_i,\g_j] \subseteq \g_{i+j}$ and the subalgebra $\g_- = \g_{-k} \dsum \cdots \dsum \g_{-1}$ is generated by $\g_{-1}$. The Lie algebra of the parabolic subgroup $P$ is the parabolic subalgebra $\p = \g_0 \dsum \cdots \dsum \g_{k}$, which has Levi decomposition $\p = \g_0 \dsum \p_+$ with $\g_0$ reductive and $\p_+ = \g_1 \dsum \cdots \dsum \g_k$ the nilradical of $\p$. At the group level, we have a reductive subgroup $G_0 \subset P$ whose Lie algebra is $\g_0$, and $P \isom G_0 \ltimes P_+$ where $P_+ = \exp (\p_+)$ is a normal, nilpotent subgroup of $P$ globally dif\/feomorphic to $\p_+$ via the exponential map.

The above-stated properties of the parabolic pair $(G,P)$ and their Lie algebras are used to identify the following important geometric structures associated to a parabolic geometry $(\Gbdle \rightarrow M,\omega)$ of type $(G,P)$. The orbit space $\Gbdle_0 := \Gbdle/P_+$ of the $P_+$-action on $\Gbdle$ def\/ines a $G_0$-principal bundle $\pi_0: \Gbdle_0 \rightarrow M$, while by def\/inition we also have a $P_+$-principal bundle $\pi_+: \Gbdle \rightarrow \Gbdle_0$. The f\/iltration of $\g$ by $\mathrm{Ad}(P)$-invariant submodules $\g^i = \g_i \dsum \cdots \dsum \g_k$ descends to a f\/iltration of $\g/\p \isom \g_-$ which is invariant under the quotient representation $\overline{\mathrm{Ad}}: P \rightarrow Gl(\g/\p)$, and thus determines a f\/iltration of the tangent bundle on the base space, $TM = T^{-k}M \supset \cdots \supset T^{-1}M$, via the isomorphism
\[
TM \isom_{\omega} \Gbdle \times_{\overline{\mathrm{Ad}}(P)} \g/\p,
\]
(which holds in general for Cartan geometries) and setting $T^iM \isom_{\omega} \Gbdle \times_{\overline{\mathrm{Ad}}(P)} \g^i/\p$. Furthermore, the Cartan connection $\omega$ descends to $\Gbdle_0$ to identify it as a reduction to $G_0$ of the structure group of the \emph{associated graded tangent bundle} $\mathrm{gr}(TM) = \mathrm{gr}_{-k}(TM) \dsum \cdots \dsum \mathrm{gr}_{-1}(TM)$, where $\mathrm{gr}_i(TM) = T^iM/T^{i+1}M$.

The data $(M,\{T^iM\},\Gbdle_0)$~-- consisting of a smooth manifold $M$, a f\/iltration $\{T^iM\}$ of its tangent bundle which satisf\/ies $\mathrm{rk}(T^iM) = \mathrm{dim}(\g^i/\p)$, and a reduction $\Gbdle_0$ of the structure group of $\mathrm{gr}(TM)$ to $G_0$ --, is called an \emph{infinitesimal flag structure of type} $(\g,P)$. The f\/lag structure is \emph{regular} if the Lie bracket of vector f\/ields on $M$ respects the f\/iltration, i.e.\ $[\Gamma(T^iM),\Gamma(T^jM)] \subset \Gamma(T^{i+j}M)$, and if the alternating bilinear form thus induced on $\mathrm{gr}(TM)$ gives it a point-wise Lie algebra structure isomorphic to $\g_-$. When the f\/lag structure is induced by a parabolic geometry of type $(G,P)$, this regularity assumption can be related to an equivalent regularity condition that the curvature of the Cartan connection $\omega$ have strictly positive homogeneity (cf.~3.1.8 of~\cite{CSbook}). A fundamental theorem of parabolic geometry states that, for any regular inf\/initesimal f\/lag structure of type $(\g,P)$, there exists a regular parabolic geometry of type $(G,P)$ which induces it. This parabolic geometry is uniquely determined up to isomorphism by a normalisation condition on the curvature of the Cartan connection, except for a number of parabolic types $(G,P)$ where the reduction of $\mathrm{gr}(TM)$ to $G_0$ provides no additional information and an extra geometric structure is needed for uniqueness (e.g.\ projective structures, where an additional choice of an equivalence class of connections is needed to f\/ix the structure). In these cases, we will assume that the extra geometric structure is included when we speak of the ``regular inf\/initesimal f\/lag structure''. We thus identify the geometric structure of an inf\/initesimal f\/lag structure with the regular, normal parabolic geometry inducing it.

In \cite{CS03}, \v{C}ap and Slov\'ak def\/ine a \emph{Weyl structure} for any parabolic geometry $(\Gbdle \rightarrow M,\omega)$ of type $(G,P)$, to be a $G_0$-equivariant section $\sigma: \Gbdle_0 \rightarrow \Gbdle$ of the $P_+$-principal bundle $\pi_+:\Gbdle \rightarrow \Gbdle_0$. We denote the set of Weyl structures by $\mathrm{Weyl}(\Gbdle,\omega)$. By Proposition~3.2 of~\cite{CS03}, global Weyl structures always exist for parabolic geometries in the real (smooth) category, and they exist locally in the holomorphic category. Considering the pull-back of the Cartan connection, $\sigma^*\omega$, the $\vert k \vert$-grading of $\g$ gives a $G_0$-invariant decomposition into components, $\sigma^*\omega = \sigma^*\omega_{-k} + \cdots + \sigma^*\omega_k,$ and by the observation that $\sigma$ commutes with fundamental vector f\/ields (i.e.\ $T_u\sigma(\widetilde{X}(u)) = \widetilde{X}(\sigma(u))$ for $X \in \g_0$ and $\widetilde{X}$ denoting the fundamental vector f\/ields of $X$ on $\Gbdle_0$ and $\Gbdle$) and the def\/ining properties of the Cartan connection, it follows that $\sigma^*\omega_i$ is horizontal for all $i \neq 0$, and that $\sigma^*\omega_0$ def\/ines a principal $G_0$ connection on $\Gbdle_0 \rightarrow M$ (cf.~3.3 of~\cite{CS03}). In particular, we see that the pair $(\Gbdle_0 \rightarrow M,\sigma^*\omega_{\leq})$ def\/ines a Cartan geometry of type $(P^*,G_0)$, where $P^* \isom \mathrm{exp}(\g_-) \rtimes G_0$ is the subgroup of $G$ containing $G_0$ with Lie algebra $\p^* = \g_- \dsum \g_0$, and the Cartan connection is given by
\[
\sigma^*\omega_{\leq} = \sigma^*\omega_{-k} + \cdots + \sigma^*\omega_0 \in \Omega^1(\Gbdle_0,\p^*).
\]

One reason Weyl structures are very useful for studying a parabolic geometry, is that they are in fact determined by very simple induced geometric objects, namely by the $\R_+$-principal connections they induce on certain ray bundles associated to $\Gbdle_0$. Fix an element $E_{\lambda}$ in the center of the reductive Lie algebra $\g_0$ such that $\mathrm{ad}(E_{\lambda})$ acts by scalar multiplication on each grading component $\g_i$ of $\g$ (for example the grading element $E$, which always exists and satisf\/ies $\mathrm{ad}(E)_{\vert \g_i} = i \cdot$). Then there is a unique representation $\lambda: G_0 \rightarrow \R_+$ satisfying $\lambda'(A) = B(E_{\lambda},A)$ for all $A \in \g_0$, $B$ the Killing form, and hence an associated $\R_+$-principal bundle $\mathcal{L}^{\lambda} \rightarrow M$. We have $\mathcal{L}^{\lambda} \isom \Gbdle_0/\mathrm{Ker}(\lambda)$ so let us denote the projection $\pi_{\lambda}: \Gbdle_0 \rightarrow \mathcal{L}^{\lambda}$. For any Weyl structure $\sigma$, the $1$-form $\lambda' \circ \sigma^*\omega_0 \in \Omega^1(\Gbdle_0)$ induces a $\R_+$-principal connection $\sigma^{\lambda}$ on $\mathcal{L}^{\lambda}$. After introducing these objects and studying their properties in Section~3 of~\cite{CS03}, \v{C}ap and Slov\'ak prove the fundamental result that the correspondence $\sigma \mapsto \sigma^{\lambda}$ def\/ines a bijective correspondence between the set of Weyl structures and the set of principal connections on $\mathcal{L}^{\lambda}$ (cf.\ Theorem~3.12 of \cite{CS03}).

{\sloppy In particular, this fact makes it possible to def\/ine certain distinguished classes of Weyl structures: A Weyl structure $\sigma$ is \emph{closed} if the induced $\R_+$-principal connection $\sigma^{\lambda}$ has vanishing curvature; it is \emph{exact} if $\sigma^{\lambda}$ is a trivial connection induced by a global trivialisation of the scale bundle $\mathcal{L}^{\lambda} \rightarrow M$. \v{C}ap and Slov\'ak prove that closed and exact Weyl structures always exist (in the smooth category), and the spaces of closed and exact Weyl structures are af\/f\/ine spaces over the closed, respectively over the exact, $1$-forms on $M$. Assuming a scale bundle $\mathcal{L}^{\lambda}$ to be f\/ixed, we denote the set of exact Weyl structures by $\overline{\mathrm{Weyl}}(\Gbdle,\omega)$, which is naturally identif\/ied with the set of global sections of $\mathcal{L}^{\lambda}$, and note that this set is non-empty (cf.\ Proposition~3.7 of~\cite{CS03}). Equivalently, an exact Weyl structure $\sigma$ is characterized by the existence of a holonomy reduction of the $G_0$-principal connection $\sigma^*\omega_0$ to the subgroup $\mathrm{Ker}(\lambda) \subset G_0$ (cf.\ Sections 3.13, 3.14 of~\cite{CS03}). We will denote this reduction by $r: \overline{\Gbdle}_0 \hookrightarrow \Gbdle_0$, and the corresponding reduction of~$\Gbdle$ to the structure group $\mathrm{Ker}(\lambda)$ by
\[
\overline{\sigma} := \sigma \circ r: \ \overline{\Gbdle}_0 \rightarrow \Gbdle.
\]
Thus an exact Weyl structure determines a Cartan geometry $(\overline{\Gbdle}_0 \rightarrow M,\overline{\sigma}^*\omega_{\leq})$ of type $(\overline{P^*},\mathrm{Ker}(\lambda))$ for $\overline{P^*} \isom \mathrm{exp}(\g_-) \rtimes \mathrm{Ker}(\lambda)$ the subgroup of $G$ containing $\mathrm{Ker}(\lambda)$ with Lie algebra $\overline{\p^*} := \g_- \dsum \mathrm{Ker}(\lambda')$. In the conformal case, exact Weyl structures correspond to metrics in the conformal equivalence class, while a general Weyl structure is given by a Weyl connection, i.e.\ a~torsion free connection which preserves the conformal equivalence class.

}

\subsection{Def\/inition and basic properties of essential structures}\label{section2.2}

Now we are ready to def\/ine essential parabolic structures and essential inf\/initesimal automorphisms. For now, let us take the following def\/initions for automorphisms, respectively inf\/initesimal automorphisms, of a Cartan geometry. For a Cartan geometry $(\Gbdle \rightarrow M,\omega)$ of arbitrary type $(G,P)$, an \emph{automorphism} $\Phi \in \mathrm{Aut}(\Gbdle,\omega)$ is a $P$-principal bundle morphism of~$\Gbdle$ such that $\Phi^*\omega = \omega$. An \emph{infinitesimal automorphism} $\mathbf{X} \in \mathrm{inf}(\Gbdle,\omega)$ is given by $\mathbf{X} \in \VF(\Gbdle)$, such that $(R_p)_*\mathbf{X} = \mathbf{X}$ and the Lie derivative satisf\/ies $\mathcal{L}_{\mathbf{X}} \omega = 0$.

Note that when $(\Gbdle,\omega)$ is a parabolic geometry, we get naturally induced $G_0$-bundle morphisms $\Phi_0:\Gbdle_0 \rightarrow \Gbdle_0$, $\R_+$-bundle morphisms $\Phi_{\lambda}: \mathcal{L}^{\lambda} \rightarrow \mathcal{L}^{\lambda}$, and dif\/feomorphisms $\varphi: M \rightarrow M$ for the $\Phi \in \mathrm{Aut}(\Gbdle,\omega)$. These are induced thanks to the $P$-equivariance of $\Phi$, by using commutativity with the appropriate projections, i.e.\ the def\/ining identities are $\Phi_0 \circ \pi_+ = \pi_+ \circ \Phi$, $\Phi_{\lambda} \circ \pi_{\lambda} = \pi_{\lambda} \circ \Phi_0$  and $\varphi \circ \pi = \pi \circ \Phi$. Hence, we get a natural action of the automorphism group $\mathrm{Aut}(\Gbdle,\omega)$ on the set of (exact) Weyl structures, by def\/ining, for $\Phi \in \mathrm{Aut}(\Gbdle,\omega)$, $\sigma \in \mathrm{Weyl}(\Gbdle,\omega)$,
\[
\Phi^*\sigma := \Phi^{-1} \circ \sigma \circ \Phi_0: \ \Gbdle_0 \rightarrow \Gbdle,
\]
and for $\sigma \in \overline{\mathrm{Weyl}}(\Gbdle,\omega)$ with corresponding global scale $s_{\sigma} \in \Gamma(\mathcal{L}^{\lambda})$,
\[
\Phi^*s_{\sigma} := \Phi_{\lambda}^{-1} \circ s_{\sigma} \circ \varphi: \ M \rightarrow \mathcal{L}^{\lambda}.
\]
Using the def\/ining relation $\pi_+ \circ \Phi^{-1} = \Phi_0^{-1} \circ \pi_+$ for $\Phi_0^{-1}$ and the corresponding relation bet\-ween~$\Phi_{\lambda}^{-1}$ and $\varphi^{-1}$, we can verify that $\Phi^*\sigma \in \Gamma(\Gbdle \rightarrow \Gbdle_0)$ and $\Phi^*s_{\sigma} \in \Gamma(\mathcal{L}^{\lambda} \rightarrow M)$; furthermore, $\Phi^*\sigma$ is $G_0$-equivariant by the equivariance properties of $\sigma$, $\Phi_0$ and $\Phi^{-1}$; i.e.\ $\Phi^*\sigma \in \mathrm{Weyl}(\Gbdle,\omega)$ and $\Phi^*s_{\sigma} \in \overline{\mathrm{Weyl}}(\Gbdle,\omega)$.

Similarly, for an inf\/initesimal automorphism $\mathbf{X} \in \VF(\Gbdle)$, we get naturally induced vector f\/ields $\mathbf{X}_0 \in \VF(\Gbdle_0)$, $\mathbf{X}_{\lambda} \in \VF(\mathcal{L}^{\lambda})$ and $X \in \VF(M)$ (with the f\/irst two being invariant with respect to the appropriate right-actions). For an arbitrary point of $M$, we may choose $\varepsilon > 0$ suf\/f\/iciently small so that $\Phi_{\mathbf{X},t}$, $\Phi_{\mathbf{X},t}^{-1}$, $\Phi_{\mathbf{X}_0,t}$, etc. all exist for $-\varepsilon \leq t \leq \varepsilon$, and so on this interval we have a~well-def\/ined family
\[
\Phi_{\mathbf{X},t}^*\sigma := \Phi_{\mathbf{X},t}^{-1} \circ \sigma \circ \Phi_{\mathbf{X}_0,t} \in \mathrm{Weyl}(\Gbdle,\omega),
\]
for any $\sigma \in \mathrm{Weyl}(\Gbdle,\omega)$, which is dif\/ferentiable in $t$ at $t=0$, so we can def\/ine the Lie derivative $\mathcal{L}_{\mathbf{X}} \sigma := (d/dt)\vert_{t=0}\Phi_{\mathbf{X},t}^*\sigma$; this is an element of the vector space on which the space of Weyl sections (an af\/f\/ine space) is modeled, i.e.\ it can be thought of as a $1$-form on $M$. For $s_{\sigma} \in \overline{\mathrm{Weyl}}(\Gbdle,\omega)$, the Lie derivative $\mathcal{L}_{\mathbf{X}}s_{\sigma}$ is def\/ined analogously, and it can be identif\/ied with an exact $1$-form on $M$.

\begin{Definition} \label{parabolic essential definition} For $(\Gbdle \rightarrow M,\omega)$ a parabolic geometry of type $(G,P)$, $\Phi \in \mathrm{Aut}(\Gbdle,\omega)$ an automorphism of the geometry, and $\sigma \in \mathrm{Weyl}(\Gbdle,\omega)$ a Weyl structure, $\Phi$ is an \emph{automorphism of~$\sigma$}, written $\Phi \in \mathrm{Aut}(\sigma)$, if and only if $\Phi^*\sigma = \sigma$ (equivalently, $\Phi \circ \sigma = \sigma \circ \Phi_0$). If $\sigma \in \overline{\mathrm{Weyl}}(\Gbdle,\omega)$ is an exact Weyl structure corresponding to $s_{\sigma} \in \Gamma(\mathcal{L}^{\lambda})$, then $\Phi$ is an \emph{exact automorphism of $\sigma$}, written $\Phi \in \mathrm{Aut}(\overline{\sigma})$, if and only if $\Phi^*s_{\sigma} = s_{\sigma}$ (equivalently, $\Phi_{\lambda} \circ s_{\sigma} = s_{\sigma} \circ \varphi$). An automorphism $\Phi \in \mathrm{Aut}(\Gbdle,\omega)$ is \emph{essential} if it is not an exact automorphism of any exact Weyl structure $\sigma \in \overline{\mathrm{Weyl}}(\Gbdle,\omega)$. We call $(\Gbdle,\omega)$ an \emph{essential parabolic structure} if $\mathrm{Aut}(\overline{\sigma}) \subsetneqq \mathrm{Aut}(\Gbdle,\omega)$ for every exact Weyl structure $\sigma$. We call a regular inf\/initesimal f\/lag structure $\mathcal{M} = (M,\{T^iM\},\Gbdle_0)$ of type $(\g,P)$ an \emph{essential structure} if the regular, normal parabolic geometry inducing it is essential.

An inf\/initesimal automorphism $\mathbf{X} \in \mathrm{inf}(\Gbdle,\omega)$ is an \emph{infinitesimal automorphism of $\sigma$}, written $\mathbf{X} \in \mathrm{inf}(\sigma)$, if and only if $\mathcal{L}_{\mathbf{X}} \sigma = 0$. For an exact Weyl structure $\sigma \in \overline{\mathrm{Weyl}}(\Gbdle,\omega)$, $\mathbf{X}$ is \emph{an exact infinitesimal automorphism of $\sigma$}, written $\mathbf{X} \in \mathrm{inf}(\overline{\sigma})$, if and only if $\mathcal{L}_{\mathbf{X}} s_{\sigma} = 0$. An inf\/initesimal automorphism is \emph{essential} if it is not an exact inf\/initesimal automorphism for any exact Weyl structure. For $\mathcal{M}$ a regular inf\/initesimal f\/lag structure as above, we say a vector f\/ield $X \in \VF(M)$ is an \emph{essential infinitesimal automorphism} of $\mathcal{M}$ if it lifts to an essential inf\/initesimal automorphism of the canonical parabolic geometry inducing $\mathcal{M}$.
\end{Definition}

\begin{Remark} Charles Frances has pointed out to us the def\/inition of essential parabolic structure given in Section~\ref{section2.2} of \cite{FrEss}, which did not make explicit use of Weyl structures. That def\/inition turns out to be almost equivalent to the one above, cf.\ Lemma \ref{inessential automorphisms TFAE}, except that in some cases the semi-simple part of the reductive group $G_0$ can be properly contained in $\mathrm{Ker}(\lambda)$. For example, in the case of strictly pseudoconvex CR structures, $G_0 \isom \R_+ \times U(n)$ and $\mathrm{Ker}(\lambda) \isom U(n)$ for an appropriate choice of scale $\lambda$. An exact Weyl structure is seen to be equivalent to a choice of pseudo-hermitian form for the CR structure, and so the above def\/inition of exact structure amounts to requiring that the group of CR transformations preserving any pseudo-hermitian form is always properly contained in the group of CR transformations.
\end{Remark}

\begin{Remark} Def\/inition \ref{parabolic essential definition} recovers the classical def\/inition of essentiality when the regular inf\/initesimal f\/lag structure is given by a conformal semi-Riemannian structure $(M,c)$ of signature $(p,q)$. In that case, $G_0$ is just the conformal group $\R_+ \times O(p,q)$, $\Gbdle_0$ is the bundle of frames which are semi-orthonormal with respect to some metric $g \in c$, and the choice of scale representation,
\begin{gather*}
 \lambda: \ \R_+ \times O(p,q) \rightarrow \R_+, \qquad
 \lambda: \ (s,A) \mapsto s^{-1},
\end{gather*}
identif\/ies $\mathcal{L}^{\lambda} \isom \Gbdle_0/\mathrm{Ker}(\lambda)$ with the ray bundle $\mathcal{Q} \rightarrow M$ of metrics in the conformal class, with the standard $\R_+$-action given by $g_x.s := s^2g_x$ for any $g \in c$ and $x \in M$ corresponding to $g_x \in \mathcal{Q}$. Exact Weyl structures thus correspond to choices of a metric in the conformal class, and a~conformal dif\/feomorphism $\varphi$ which uniquely corresponds to an automorphism $\Phi \in \mathrm{Aut}(\Gbdle,\omega)$ of the canonical conformal Cartan geometry is essential in the sense of Def\/inition~\ref{parabolic essential definition} if $\Phi \notin \mathrm{Aut}(\overline{\sigma})$ for all exact Weyl structures $\sigma$, i.e.\ if $\varphi$ fails to preserve all metrics in the conformal class.
\end{Remark}

\begin{Remark} \label{underlying automorphisms} Let $\mathcal{M} = (M,\{T^iM\},\Gbdle_0)$ denote a regular inf\/initesimal f\/lag structure of some parabolic type $(\g,P)$, if necessary including the extra geometric data required so that the regular, normal parabolic geometry of type $(G,P)$ inducing it is unique up to isomorphism. If we need to distinguish this parabolic geometry from others of the same type, we will use the notation $(\Gbdle,\omega^{nc})$ to signify the canonical (\textbf{n}ormal \textbf{C}artan) geometry. In this setting, we can def\/ine an automorphism of the structure in terms of $\mathcal{M}$: An automorphism of the regular inf\/initesimal f\/lag structure, $\varphi \in \mathrm{Aut}(\mathcal{M})$, is a dif\/feomorphism $\varphi \in \mathrm{Dif\/f}(M)$ which satisf\/ies: $(i)$~$\varphi_*(T_x^iM) \subseteq T_{\varphi(x)}^iM$ for all $x \in M$ and all $-k \leq i \leq -1$; and $(ii)$~the induced bundle map $\mathrm{gr}(\varphi)$ (which as a~consequence of $(i)$~is a lift of $\varphi$ def\/ined on the bundle $\mathcal{F}(\mathrm{gr}(TM))$ of frames of the associated graded tangent bundle) preserves $\Gbdle_0$ as a subbundle of $\mathcal{F}(\mathrm{gr}(TM))$ (and hence $\mathrm{gr}(\varphi)$ restricts to a $G_0$-bundle morphism $\Phi_0$ of $\Gbdle_0$). We can identify $\mathrm{Aut}(\mathcal{M})$ with $\mathrm{Aut}(\Gbdle,\omega^{nc})$ since by uniqueness of $(\Gbdle,\omega^{nc})$ up to isomorphism, $\varphi$ (and $\Phi_0$) lift to a unique $P$-bundle morphism $\Phi$ of $\Gbdle$ preserving $\omega^{nc}$ under pullback. Thus, we can think of an automorphism $\varphi \in \mathrm{Aut}(\mathcal{M})$ as including as well the automorphism $\Phi \in \mathrm{Aut}(\Gbdle,\omega^{nc})$ and the induced $G_0$-bundle morphism $\Phi_0 = \mathrm{gr}(\varphi)_{\vert \Gbdle_0}$ of $\Gbdle_0$.
\end{Remark}

\begin{Lemma} \label{inessential automorphisms TFAE} Let $\Phi \in \mathrm{Aut}(\Gbdle,\omega)$ be an automorphism of a parabolic geometry, let $\sigma$ be a Weyl structure, and let a scale bundle $\mathcal{L}^{\lambda} \rightarrow M$ be fixed. The following are equivalent:
\begin{enumerate}\itemsep=0pt
\item[$(i)$] $\Phi \in \mathrm{Aut}(\sigma)$;
\item[$(ii)$] For the induced bundle morphism $\Phi_0: \Gbdle_0 \rightarrow \Gbdle_0$, we have $\Phi_0 \in \mathrm{Aut}(\Gbdle_0,\sigma^*\omega_{\leq})$;
\item[$(iii)$] $\Phi_{\lambda}$ preserves the scale bundle connection $\sigma^{\lambda} \in \Omega^1(\mathcal{L}^{\lambda})$: $\Phi_{\lambda}^*\sigma^{\lambda} = \sigma^{\lambda}$.
\end{enumerate}

 If $\sigma$ is exact, then $\Phi \in \mathrm{Aut}(\overline{\sigma})$ if and only if the induced bundle morphism $\Phi_0$ preserves the sub-bundle $\overline{\Gbdle}_0 \subset \Gbdle_0$ and the restriction satisfies $(\Phi_0)_{\vert \overline{\Gbdle}_0} \in \mathrm{Aut}(\overline{\Gbdle}_0,\overline{\sigma}^*\omega_{\leq})$.

For $\mathbf{X} \in \mathrm{inf}(\Gbdle,\omega)$, $\mathbf{X} \in \mathrm{inf}(\sigma)$ if and only if $\mathbf{X}_0 \in \mathrm{inf}(\Gbdle_0,\sigma^*\omega_{\leq})$ and, for $\sigma$ exact, $\mathbf{X} \in \mathrm{inf}(\overline{\sigma})$ if and only if the restriction of $\mathbf{X}_0$ to $\overline{\Gbdle}_0$ is tangent to $\overline{\Gbdle}_0$ and induces an element of $\mathrm{inf}(\overline{\Gbdle}_0,\overline{\sigma}^*\omega_{\leq})$.
\end{Lemma}

\begin{proof} $(i)$ $\Rightarrow$ $(ii)$, since by $\Phi \circ \sigma = \sigma \circ \Phi_0$ and $\Phi \in \mathrm{Aut}(\Gbdle,\omega)$ we have $\Phi_0^*(\sigma^*\omega) = \sigma^*\omega$, and in particular $\Phi_0^*(\sigma^*\omega_{\leq}) = \sigma^*\omega_{\leq}$. And $(ii)$ $\Rightarrow$ $(iii)$, since $\Phi_0^*(\sigma^*\omega_{\leq}) = \sigma^*\omega_{\leq}$ implies that $\Phi_0^*(\sigma^*\omega_0) = \sigma^*\omega_0$. In particular, $\Phi_0^*(\lambda' \circ \sigma^*\omega_0) = \lambda' \circ \Phi_0^*(\sigma^*\omega_0) = \lambda' \circ \sigma^*\omega_0$. Hence, the $\R_+$-bundle morphism $\Phi_{\lambda}$ and the $\R_+$-principal connection $\sigma^{\lambda}$, induced on $\mathcal{L}^{\lambda}$ by $\Phi_0$ and $\lambda' \circ \sigma^*\omega_0$, respectively, satisfy: $\Phi_{\lambda}^* \sigma^{\lambda} = \sigma^{\lambda}$.

  We now show $(iii)$ $\Rightarrow$ $(i)$: Consider the $\R_+$-principal connection $(\Phi^*\sigma)^{\lambda} \in \Omega^1(\mathcal{L}^{\lambda})$. It is induced by:
\begin{gather*}
\lambda' \circ (\Phi^*\sigma)^*\omega_0  = \lambda' \circ \big(\big(\Phi^{-1}\circ \sigma \circ \Phi_0\big)^*\omega_0\big)
  = \lambda' \circ \big(\Phi_0^*\sigma^*\big(\Phi^{-1}\big)^*\omega_0\big) \\
\phantom{\lambda' \circ (\Phi^*\sigma)^*\omega_0}{} = \lambda' \circ (\Phi_0^*\sigma^*\omega_0) = \Phi_0^*(\lambda' \circ \sigma^*\omega_0).
\end{gather*}
So $(\Phi^*\sigma)^{\lambda} = \Phi_{\lambda}^*\sigma^{\lambda}$, which equals $\sigma^{\lambda}$ by assumption $(iii)$. Thus, by Theorem~3.12 of \cite{CS03} (cf.\ discussion in Section~\ref{section2.1}), the Weyl structures $\Phi^*\sigma$ and $\sigma$ are equal, showing that $(i)$ holds.

 To see the f\/inal statement of the lemma, let $\sigma$ be an exact Weyl structure and let us denote by $s_{\sigma} \in \Gamma(\mathcal{L}^{\lambda})$ the global scale which induces the trivial connection $\sigma^{\lambda} \in \Omega^1(\mathcal{L}^{\lambda})$. That is, for any point $p = s_{\sigma}(x).r \in \mathcal{L}^{\lambda}$, for $x \in M$, $r \in \R_+$, we have the decomposition
\[
T_p \mathcal{L}^{\lambda} = (R_r)_*((s_{\sigma})_*(T_xM)) \dsum \R \zeta_1(p),
\]
for $\zeta_1$ the fundamental vector f\/ield on $\mathcal{L}^{\lambda}$ of the vector $1 \in \R$; the value of $\sigma^{\lambda}$ on a tangent vector $v \in T_p\mathcal{L}^{\lambda}$ is given by the coef\/f\/icient of $\zeta_1(p)$ determined by this decomposition. Then the holonomy reduction of $(\mathcal{L}^{\lambda},\sigma^{\lambda})$ to the trivial structure group is given by $s_{\sigma}(M) \subset \mathcal{L}^{\lambda}$, and the reduction of $(\Gbdle_0,\sigma^*\omega_0)$ to $\mathrm{Ker}(\lambda)$ is given by
\[
\overline{\Gbdle}_0 = \pi_{\lambda}^{-1}(s_{\sigma}(M)) \subset \Gbdle_0.
\]
Hence for $x \in M$ and $u \in (\overline{\Gbdle}_0)_x$, we have $\pi_{\lambda}(u) = s_{\sigma}(x)$ and $\Phi_0(u) \in \overline{\Gbdle}_0$ if and only if $\pi_{\lambda}(\Phi_0(u)) = s_{\sigma}(\pi_0(\Phi_0(u)))$. But $\pi_{\lambda}(\Phi_0(u)) = \Phi_{\lambda}(\pi_{\lambda}(u)) = \Phi_{\lambda}(s_{\sigma}(x))$ and $\pi_0(\Phi_0(u)) = \varphi(\pi_0(u)) = \varphi(x)$. Thus, $\Phi_0(u) \in \overline{\Gbdle}_0$ if and only if $\Phi_{\lambda}(s_{\sigma}(x)) = s_{\sigma}(\varphi(x))$. But if $\Phi \in \mathrm{Aut}(\overline{\sigma})$, then clearly $\Phi_{\lambda}$ preserves the induced (trivial) connection $\sigma^{\lambda} \in \Omega^1(\mathcal{L}^{\lambda})$, so by $(iii)$ $\Leftrightarrow$ $(ii)$ shown above we have $\Phi_0^*(\sigma^*\omega_{\leq}) = \sigma^*\omega_{\leq}$ and since $\Phi_0$ preserves $\overline{\Gbdle}_0$ the corresponding identity follows for the restriction and $\overline{\sigma}^*\omega_{\leq}$.

 The statements for $\mathbf{X} \in \mathrm{inf}(\Gbdle,\omega)$ are proven in the same manner.
 \end{proof}

\begin{Remark} In particular, it follows from the proof of Lemma \ref{inessential automorphisms TFAE} that any exact automor\-phism~$\Phi$ of an exact Weyl structure $\sigma$ is also an automorphism of $\sigma$, i.e.\ $\Phi \in \mathrm{Aut}(\overline{\sigma}) \Rightarrow \Phi \in \mathrm{Aut}(\sigma)$ as one would hope.

 On the other hand, the converse does not hold: If $\Phi \in \mathrm{Aut}(\sigma)$ for an exact Weyl structure $\sigma$, the requirement that $\Phi_0(\overline{\Gbdle}_0) \subset \overline{\Gbdle}_0$ is necessary to guarantee that an automorphism $\Phi$ of an exact Weyl structure $\sigma$ is in fact an \emph{exact} automorphism. An instructive example is the conformal structure induced by the Euclidean metric on $\R^n$: The dif\/feomorphism given by dilation by a~positive constant $r$ is always an automorphism of the exact Weyl structure corresponding to the Euclidean metric. However, for $r \neq 1$, this dif\/feomorphism is not an isometry, hence not an \emph{exact} automorphism. We are grateful to Felipe Leitner for bringing this to our attention, which led us to modify an earlier version of Def\/inition~\ref{parabolic essential definition}.
\end{Remark}

\section{Lichn\'{e}rowicz theorem for rank one parabolic geometries}\label{section3}

In this section, we establish Theorem~\ref{Lich for rank one parabolic} for the so-called rank one parabolic geometries. These are parabolic geometries of types~$(G,P)$ such that the homogeneous model~$G/P$ is the bounda\-ry of a real rank one symmetric space~$G/K$ (for $K$ a maximal compact subgroup). The boundary~$G/P$ is dif\/feomorphic to the sphere of dimension $\mathrm{dim}(G/K)-1$, while the symmetric spaces~$G/K$ are given by the hyperbolic spaces of real, complex and quaternionic type in all appropriate real dimensions, and the hyperbolic Cayley plane in real dimension $16$. As mentioned earlier, the rank one parabolic geometries associated to the corresponding types are, respectively, conformal Riemannian structures, strictly pseudo-convex partially-integrable CR structures, quaternionic contact structures (of positive-def\/inite signature), and octonionic contact structures. The key result needed to prove this theorem is the following, proved by C.~Frances (Theorem~3 of~\cite{Frances}), which generalizes theorems of Ferrand~\cite{Ferrand} and Schoen~\cite{Schoen} in the cases of conformal Riemannian and strictly pseudo-convex CR structures:

\begin{Theorem}[Frances, \cite{Frances}]  \label{Frances thm}Let $(\Gbdle \rightarrow M,\omega)$ be a regular rank one parabolic geometry, with~$M$ connected. If $\mathrm{Aut}(\Gbdle,\omega)$ acts improperly on $M$, then $M$ is geometrically isomorphic to either the compact homogeneous model $G/P$ or the noncompact space $G/P \backslash \{eP\}$.
\end{Theorem}

Here, ``geometrically isomorphic'' means there is a dif\/feomorphism of $M$ onto the space in question, which is covered by a morphism of Cartan bundles which pulls back the Maurer--Cartan connection to~$\omega$. The assumption that the parabolic geometries are of rank one type is key for the argument in~\cite{Frances}, because it allows the author to exploit the so-called ``north-south dynamics'' on the homogeneous model $G/P$. Theorem~\ref{Lich for rank one parabolic} now follows as a result of Theorem~\ref{Frances thm} and the following proposition, which generalizes a result of~\cite{Aleks} and whose proof follows the same line of argumentation:

\begin{Proposition} \label{proper implies inessential} If $(\Gbdle \rightarrow M,\omega)$ is an essential parabolic structure, then $\mathrm{Aut}(\Gbdle,\omega)$ acts improperly on~$M$.
\end{Proposition}

\begin{proof} Fix a bundle of scales $\mathcal{L}^{\lambda} \rightarrow M$ for $(\Gbdle,\omega)$. Assume that $\mathrm{Aut}(\Gbdle,\omega)$ acts properly on~$M$ and let us show that the parabolic structure is not essential. By def\/inition, it suf\/f\/ices to construct a global scale $s: M \rightarrow \mathcal{L}^{\lambda}$ which is $\mathrm{Aut}(\Gbdle,\omega)$-invariant, i.e.\ such that $\Phi_{\lambda} \circ s = s \circ \varphi$ holds for all $\Phi \in \mathrm{Aut}(\Gbdle,\omega)$ and $\Phi_{\lambda}: \mathcal{L}^{\lambda} \rightarrow \mathcal{L}^{\lambda}$, $\varphi: M \rightarrow M$ the induced dif\/feomorphisms.

 We construct this $\mathrm{Aut}(\Gbdle,\omega)$-invariant scale $s$ using classical properties of proper group actions. The so called ``tube theorem'' (alias ``slice theorem'', cf.\ e.g.\ Theorem~2.4.1 in~\cite{DuistKolk}) guarantees the following, for a $C^{\infty}$-action of a Lie group $H$ on a manifold $M$ which is proper at $x \in M$: There exists a $H$-invariant neighborhood $U$ of $x$ on which the $H$-action is equivalent to the left $H$-action on the quotient space $H \times_{K} B$~-- for $K \subset H$ a compact subgroup and $B$ a $K$-invariant neighborhood of $0$ in a $K$-module~$V$~-- given by $h_1.[h_2,b] = [h_1h_2,b]$ for $h_i \in H, b \in B$ and $[h,b]$ the equivalence class of $(h,b) \in H \times B$ under the left $K$-action $k.(h,b) := (h.k^{-1},k.b)$. Starting from a choice of global scale $s_0: M \rightarrow \mathcal{L}^{\lambda}$, and letting $H = \mathrm{Aut}(\Gbdle,\omega)$, $e \in H$ the identity automorphism and $\Phi \in H$ arbitrary, set:
\begin{gather}
s_{U}([e,b]) := \int_{\Psi \in K} (\Psi_{\lambda})^{-1}(s_0([\Psi,b]))d\Psi; \label{star}\\
s_U([\Phi,b]) := \Phi_{\lambda}(s_U([e,b])). \label{star star}
\end{gather}

 One verif\/ies that this gives a well-def\/ined local section $s_{U}: U \rightarrow \mathcal{L}^{\lambda}_{\vert U}$, which involves checking that for $(e,b) \sim (\Phi,b')$ (i.e.\ for $\Phi \in K$ and $b = \varphi(b')$) the values $s_U([e,b])$ given by~(\ref{star}) and $s_U([\Phi,b'])$ given by~(\ref{star star}), agree. This follows by unwinding the def\/initions, and using a bi-invariant Haar measure~$d\Psi$ on the compact group~$K$. And since $[\Phi,b] = \Phi.[e,b]$ corresponds to the point $\varphi(x')$ for $x' \simeq [e,b]$, the def\/ining equation~(\ref{star star}) automatically gives us the invariance property, $s_U \circ \varphi = \Phi_{\lambda} \circ s_U$.

 We go from local $\mathrm{Aut}(\Gbdle,\omega)$-invariant sections $s_U: U \rightarrow \mathcal{L}^{\lambda}_U$ to a global $\mathrm{Aut}(\Gbdle,\omega)$-invariant scale $s_{\rm inv} \in \Gamma(\mathcal{L}^{\lambda})$ as follows: First, note that $\mathcal{L}^{\lambda}$ admits a global section $s \in \Gamma(\mathcal{L}^{\lambda})$ and hence a global trivialisation $\mathcal{L}^{\lambda} \isom_s M \times \R_+$ by identifying $s(x).r \simeq_s (x,r)$ for any $x \in M$ and $r \in \R_+$. This in turn induces a bijection $r_s: \Gamma(\mathcal{L}^{\lambda}) \rightarrow C^{\infty}(M,\R_+)$. For any $\Phi \in \mathrm{Aut}(\Gbdle,\omega)$, def\/ine the smooth function $\varphi_s \in C^{\infty}(M,\R_+)$ by the identity $\Phi_{\lambda}(s(x)) = s(\varphi(x)).\varphi_s(x)$. Then we see that a $\Phi$-invariant scale $s' \in \Gamma(\mathcal{L}^{\lambda})$ corresponds, under the bijection $r_s$, to a smooth function $r' = r_s(s') \in C^{\infty}(M,\R_+)$ which satisf\/ies $r'(\varphi(x)) = \varphi_s(x)r'(x)$.

 Next, note that since $\mathrm{Aut}(\Gbdle,\omega)$ acts properly on $M$, there exists a covering $\{U_{\alpha}\}$ of $M$ by $\mathrm{Aut}(\Gbdle,\omega)$-invariant open sets as above (so they all admit invariant local scales $s_{\alpha}:U_{\alpha} \rightarrow \mathcal{L}^{\lambda}_{U_{\alpha}}$ and we denote $r_{\alpha} := r_s(s_{\alpha}) \in C^{\infty}(U_{\alpha},\R_+)$); and (applying Theorem 6 of \cite{Aleks}) there exists a~partition of unity $\{f_i\}$ for $M$ which is subordinate to $\{U_{\alpha}\}$ (so $\mathrm{supp}(f_i) \subset U_{\alpha(i)}$ for all $i$) and the $f_i$ are $\mathrm{Aut}(\Gbdle,\omega)$-invariant. Now we def\/ine $s_{\rm inv} \in \Gamma(\mathcal{L}^{\lambda})$ via $r_{\rm inv} = r_s(s_{\rm inv}) \in C^{\infty}(M,\R_+)$ and the formula
\[
r_{\rm inv}(x) = \sum_{i}f_i(x)r_{\alpha(i)}(x)
\]
for all $x \in M$. This is well-def\/ined, smooth, and from the $\mathrm{Aut}(\Gbdle,\omega)$-invariance of the $f_i$, together with the $\mathrm{Aut}(\Gbdle,\omega)$-invariance of the local sections $s_{\alpha}$, we compute that $r_{\rm inv}(\varphi(x)) = \varphi_s(x)r_{\rm inv}(x)$ for all $x \in M$ and all $\Phi \in \mathrm{Aut}(\Gbdle,\omega)$, i.e.\ $s_{\rm inv} \in \Gamma(\mathcal{L}^{\lambda})$ is a $\mathrm{Aut}(\Gbdle,\omega)$-invariant global scale.
\end{proof}

\section{Proof of local results}\label{section4}

A key reference for the study of inf\/initesimal automorphisms of parabolic geometries is \cite{Cap inf aut}. In that text, A.~\v{C}ap generalized to arbitrary parabolic geometries a bijective correspondence between conformal vector f\/ields and adjoint tractors (sections of the associated bundle to the canonical Cartan bundle, $\Gbdle \rightarrow M$, induced by the adjoint representation on $\g$) satisfying an identity involving the Cartan curvature, which was f\/irst discovered by A.R.~Gover in~\cite{Gover}. Moreover, the text of \v{C}ap relates this general bijective correspondence to the f\/irst splitting operator of a so-called curved BGG-sequence for the parabolic geometry, cf.\ Theorem~3.4 of~\cite{Cap inf aut}. The curvature identity of \cite{Cap inf aut} extends without dif\/f\/iculty to general inf\/initesimal automorphisms of Cartan geometries. This allows us to apply this fundamental identity to the Cartan geometries $(\Gbdle_0,\sigma^*\omega_{\leq})$, which we do in Section~\ref{section4.2} to establish a general ``dictionary'' between essentiality of an inf\/initesimal automorphism near a singularity, and the so-called holonomy associated to such a singularity (cf.\ Def\/inition~\ref{holonomy definition}; again, this should not be confused with the holonomy of the Cartan connection).

\subsection{Background results on inf\/initesimal automorphisms}\label{section4.1}

We recall some general notions, mainly following the development of \cite{Cap inf aut} (cf.\ also~1.5 of \cite{CSbook}), in the setting of a general Cartan geometry $(\Gbdle \rightarrow M,\omega)$ of type $(G,P)$ (for now not assumed to be of parabolic type). For any representation $\rho: P \rightarrow Gl(V)$, we have the associated vector bundle $V(M) := \Gbdle \times_{\rho} V$. The smooth sections of such a bundle are identif\/ied with $P$-equivariant, $V$-valued smooth functions on $\Gbdle$ in the standard manner, and we will simply treat them as such:
\[
\Gamma(V(M)) = \big\{f \in C^{\infty}(\Gbdle,V) \, \vert \, f(u.p) = \rho(p^{-1})(f(u))\big\} =: C^{\infty}(\Gbdle,V)^P.
\]
For the most part, the important associated bundles we are dealing with are \emph{tractor bundles}, which for our purposes simply means that the representation $(\rho,V)$ is the restriction to $P$ of a $G$-representation $\tilde{\rho}:G \rightarrow Gl(V)$. And the primary tractor bundle is the adjoint bundle induced by the restriction of the adjoint representation $\mathrm{Ad}: G \rightarrow Gl(\g)$ to $P$, which we will denote by $\Abdle = \Abdle(M)$ if there is no danger of confusion about which Lie algebra $\g$ is meant, and otherwise by $\g(M)$. Note that the Lie bracket $[\cdot \, ,\cdot]_{\g}$ of $\g$, by $\mathrm{Ad}(P)$-invariance, determines an algebraic bracket on f\/ibers of $\Abdle$ as well as on sections, which we denote with curly brackets $\{\cdot \, ,\cdot\}: \Abdle \times \Abdle \rightarrow \Abdle$. Also, note that there is a natural projection $\Pi: \Abdle \rightarrow TM$, induced by the projection $\g \rightarrow \g/\p$ and the isomorphism $TM \isom_{\omega} \Gbdle \times_{\overline{\mathrm{Ad}}(P)} \g/\p$.

The Cartan connection determines an identif\/ication of right-invariant vector f\/ields
\[
\VF(\Gbdle)^P = \{ \mathbf{X} \in \VF(\Gbdle) \, \vert \, \mathbf{X}(u.p) = (R_p)_* (\mathbf{X}(u))\},
\]
with sections of the adjoint bundle. Namely, to $\mathbf{X} \in \VF(\Gbdle)$ we associate a function $s_{\mathbf{X}} \in C^{\infty}(\Gbdle,\g)$ def\/ined by $s_{\mathbf{X}}(u) := \omega(\mathbf{X}(u))$; conversely, to a function $s \in C^{\infty}(\Gbdle,\g)$, associate $\mathbf{X}_s \in \VF(\Gbdle)$ def\/ined by $\mathbf{X}_s(u) = \omega_u^{-1}(s(u))$. The property~(\ref{parallelism}) of a~Cartan connection insures that both maps are well-def\/ined, they are inverse, and by property~(\ref{ad-inv}) of $\omega$ these maps restrict to an isomorphism $\VF(\Gbdle)^P \isom_{\omega} \Gamma(\Abdle)$. Similarly, the Cartan connection $\omega$ induces natural identif\/ications $\overline{\Omega^k(\Gbdle;\g)}^P \isom_{\omega} \Omega^k(M;\Abdle)$ of the horizontal, $\mathrm{Ad}(P)$-equivariant $\g$-valued $k$-forms on $\Gbdle$ with the $\Abdle$-valued $k$-forms on~$M$. For a tractor bundle~$V(M)$, the identif\/ication $\VF(\Gbdle)^P \isom_{\omega} \Gamma(\Abdle)$ yields two kinds of dif\/ferentiation of smooth sections with respect to adjoint tractors:

\begin{Definition} \label{differentiation} The \textit{invariant differentiation} or \textit{fundamental D-operator} of $V(M)$ is the map $D^{V} : \Gamma(V(M)) \rightarrow \Gamma(\Abdle^* \tens V(M))$ def\/ined, for any $s \in \Gamma(\Abdle)$ and any $v \in \Gamma(V(M))$, by:
\[
D_s^{V}v := \mathbf{X}_s(v).
\]
The \textit{tractor connection} of $V(M)$ is the map $\nabla^{V}: \Gamma(V(M)) \rightarrow \Gamma(\Abdle^* \tens V(M))$ def\/ined, for any $s \in \Gamma(\Abdle)$ and any $v \in \Gamma(V(M))$, by:
\begin{gather}
\nabla^{V}_s v := D^{V}_s v + (d\tilde{\rho} \circ s) \circ v. \label{tractor connection}
\end{gather}
 Recall that $\tilde{\rho}: G \rightarrow Gl(V)$ is the $G$-representation which~$\rho$ is a restriction of, given by the def\/inition of a tractor bundle. In fact, the quantity def\/ined by (\ref{tractor connection}) only depends on the equiva\-len\-ce class~$[s]$ of~$s$ under the quotient $\Abdle/\p(M) = \Gbdle \times_{\overline{\mathrm{Ad}}(P)} \g/\p \isom TM$, and we identify the tractor connection with a covariant derivative on~$V(M)$:
\[
\nabla^V: \ \Gamma(V(M)) \rightarrow \Gamma(T^*M \tens V(M)).
\]
\end{Definition}

The curvature tensor of a Cartan connection is the $\g$-valued two-form on $\Gbdle$ def\/ined, for any $u \in \Gbdle$ and $v, w \in T_u\Gbdle$, by the structure equation $\Omega^{\omega}(v,w) := d\omega(v,w) + [\omega(v),\omega(w)]$. The curvature tensor is horizontal and $P$-equivariant, and we may equivalently consider the curvature function $\kappa^{\omega} \in C^{\infty}(\Gbdle;\Lambda^2(\g/\p)^* \tens \g)^P$ induced by $\kappa^{\omega}(u)(X,Y) := \Omega^{\omega}(\omega_u^{-1}(X),\omega_u^{-1}(Y))$ for any $u \in \Gbdle$ and $X,Y \in \g$. The following identity was proven for parabolic geometries in~\cite{Cap inf aut} and the proof carries over without substantive changes to Cartan geometries of general type (cf.\ also Lemma~1.5.12 in~\cite{CSbook}):

\begin{Lemma} \label{infinitesimal deformation} Let $\mathbf{X} \simeq_{\omega} s_{\mathbf{X}}$ for $\mathbf{X} \in \VF(\Gbdle)^P$ and $s_{\mathbf{X}} \in \Gamma(\Abdle)$. Then for $\mathcal{L}_{\mathbf{X}} \omega \in \overline{\Omega^1(\Gbdle;\g)}^P$ we have, under the identification $\overline{\Omega^1(\Gbdle;\g)}^P \isom_{\omega} \Omega^1(M;\Abdle)$:
\begin{gather}
\mathcal{L}_{\mathbf{X}} \omega \simeq_{\omega} \nabla^{\Abdle} s_{\mathbf{X}} + \Pi(s_{\mathbf{X}}) \: \k^{\omega}. \label{infinitesimal deformation identity}
\end{gather}
\end{Lemma}

\subsection{Holonomy and essentiality of inf\/initesimal automorphisms}\label{section4.2}

Let us return now to the setting of a (regular, normal) parabolic geometry $(\Gbdle \rightarrow M,\omega)$ of type $(G,P)$ and the corresponding regular inf\/initesimal f\/lag structure $\mathcal{M} = (M,\{T^iM\},\Gbdle_0)$ of type $(\g,P)$. As was noted for automorphisms of $(\Gbdle,\omega^{nc})$ in Remark \ref{underlying automorphisms}, we may determine an inf\/initesimal automorphism $\mathbf{X} \in \mathrm{inf}(\Gbdle,\omega)$ by conditions on the underlying vector f\/ield $X \in \VF(M)$, which just amount to imposing the same conditions for the locally def\/ined dif\/feomorphisms given by f\/lowing along $X$, i.e.\ we must have $[X,\Gamma(T^iM)] \subseteq \Gamma(T^iM)$ for all $-k \leq i \leq -1$, and the condition that the local f\/lows of $X$ determine local bundle maps of $\mathcal{F}(\mathrm{gr}(TM))$ which preserve~$\Gbdle_0$ as a~subbundle. In the dif\/ferent examples of parabolic geometries, this translates into more geometric language. For example, in the conformal case, the former condition is trivial, while the latter condition amounts to requiring the conformal Killing equation, $\Lie_X g = \lambda g$ for any $g \in c$ and some $\lambda = \lambda(X) \in C^{\infty}(M)$. In the case of CR structures, the conditions are that $[X,\Gamma(\mathcal{H})] \subseteq \Gamma(\mathcal{H})$ for $\mathcal{H} \subset TM$ the codimension one contact distribution def\/ining the CR structure, and that $\Lie_X J = 0$ for $J$ the almost complex structure on $\mathcal{H}$. We write $X \in \mathrm{inf}(\mathcal{M})$ and consider the lift $\mathbf{X} \in \mathrm{inf}(\Gbdle,\omega^{nc})$ to be implicitly included.

Using Lemma~\ref{infinitesimal deformation} and Lemma~\ref{inessential automorphisms TFAE}, we obtain a bijection between inf\/initesimal automorphisms $\mathbf{X} \in \mathrm{inf}(\Gbdle,\omega)$ and adjoint tractors $s_{\mathbf{X}} \in \Gamma(\Abdle)$ satisfying $\nabla^{\Abdle}s_{\mathbf{X}} + \Pi(s_{\mathbf{X}}) \: \k^{\omega} = 0$. In the present setting, this gives us a bijection between $X \in \mathrm{inf}(\mathcal{M})$ and such $s_{\mathbf{X}}$, and moreover it is easy to verify that $\Pi(s_{\mathbf{X}}) = X \in \VF(M)$. Now, denote by $\mathbf{X}_0 \in \VF(\Gbdle_0)^{G_0}$ the invariant vector f\/ield induced, via projection by $\pi_+$, by $\mathbf{X} \in \VF(\Gbdle)^P$. For any Weyl structure $\sigma: \Gbdle_0 \rightarrow \Gbdle$, the induced Cartan connection $\sigma^*\omega_{\leq}$ gives us an isomorphism denoted $\VF(\Gbdle_0)^{G_0} \isom_{\sigma} \Gamma(\Abdle^{\sigma})$, where we denote with $\Abdle^{\sigma} = \p^*(M)$ the adjoint bundle of the Cartan geometry $(\Gbdle_0 \rightarrow M,\sigma^*\omega_{\leq})$. Let us write $\mathbf{X}_0 \simeq_{\sigma} s_{\mathbf{X}_0}$ for the adjoint tractor corresponding to $\mathbf{X}_0 \in \VF(\Gbdle_0)^{G_0}$. If we further denote by $\nabla^{\sigma}: \Gamma(\Abdle^{\sigma}) \rightarrow \Gamma(T^*M\tens\Abdle^{\sigma})$ the corresponding adjoint tractor connection, then Lemma \ref{infinitesimal deformation} tells us that $X \in \mathrm{inf}(\mathcal{M})$ is an inf\/initesimal automorphism of some Weyl structure $\sigma: \Gbdle_0 \rightarrow \Gbdle$ if and only if
\[
\nabla^{\sigma} s_{\mathbf{X}_0} + X \: \k^{\sigma^*\omega_{\leq}} = 0.
\]

At present, we are only interested in local properties of inf\/initesimal automorphisms, viz the question if some neighborhood of a given point can be found on which $X$ is inessential. The following proposition shows that the points $x \in M$ for which the answer could be ``no'', must all be singularities of the vector f\/ield, i.e.\ $X(x)=0$ (where, incidentally, the above identity simplif\/ies to $(\nabla^{\sigma} s_{\mathbf{X}_0})(x) = 0$).

\begin{Proposition} \label{non-singular points} Let $\mathcal{M} = (M,\{T^iM\},\Gbdle_0)$ be a regular infinitesimal flag structure of type $(\g,P)$, let $X \in \mathrm{inf}(\mathcal{M})$, and let $x \in M$. If $X(x) \neq 0$, then there exists a neighborhood $U$ of $x$ such that the restriction of $X$ to $U$ is inessential.
\end{Proposition}

\begin{proof} Take a neighborhood $U$ of $x$ on which f\/low-box coordinates for the f\/low of $X$ can be introduced:
\[
U = \{(x_0,\varphi_{X,t}(x_0)) \, \vert \, x_0 \in M_0 \cap U, \, -\varepsilon < t < \varepsilon\},
\]
where $M_0$ is some locally def\/ined hypersurface transversal to the integral curves $\varphi_{X,t}(x_0)$ of~$X$, which are def\/ined for the interval given. This can be done since we may f\/irst restrict an open neighborhood of $x$ on which $X$ is non-vanishing. But this provides all the features needed to transfer the argument used in the proof of Proposition \ref{proper implies inessential} to establish the existence of local $\mathrm{Aut}(\Gbdle,\omega)$-invariant sections $s_U: U \rightarrow \mathcal{L}^{\lambda}_{\vert U}$ to the current context: Namely, the formula $s(x_0,\varphi_{X,t}(x_0)) :=  \Phi_{\mathbf{X}_{\lambda},t}(s(x_0))$ gives a well-def\/ined local scale $s: U \rightarrow \mathcal{L}^{\lambda}_{\vert U}$ which by construction is invariant under the f\/lows $\varphi_{X,t}$ for $t$ suf\/f\/iciently small. \end{proof}

From now on, let us f\/ix a singularity $x \in M$ of an inf\/initesimal automorphism $X \in \mathrm{inf}(\mathcal{M})$. We also choose a point $u \in \Gbdle$ in the f\/iber over $x$, and let $u_0 = \pi_+(u) \in \Gbdle_0$, likewise in the f\/iber over $x$. The remaining text is aimed at relating the local essentiality of $X$ near $x$, to invariant properties of the holonomy of $X$ at $x$, which is a one-parameter subgroup $h^t \subset P$:

\begin{Definition}[cf.~\cite{Frances3}, Section~6] \label{holonomy definition} Given $X$, $x$ and $u$ as above, \emph{the holonomy $h_u^t$ of $X$ at $x$ with respect to $u$} is def\/ined, for $t$ suf\/f\/iciently small, as follows: Let $\Phi_{\mathbf{X},t}(u)$, the integral curve of $\mathbf{X}$ through $u$, be def\/ined for $t \in (-\varepsilon,\varepsilon)$. Since $\mathbf{X}$ projects to $X$, $\mathbf{X}(u')$ is tangent to $\Gbdle_x$ for all $u' \in \Gbdle_x$, and hence all $\Phi_{\mathbf{X},t}(u)$ lie in $\Gbdle_x$. Then $h_u^t \in P$ is def\/ined by:
\[
\Phi_{\mathbf{X},t}(u) =: u.h^t_u.
\]
Since $h_u^{t+s} = h_u^th_u^s$ whenever both are def\/ined, $h_u^t = \mathrm{exp}(tX_{h,u})$ for some $X_{h,u} \in \p$, and we def\/ine~$h_u^t$ via this identity for all $t \in \R$.
\end{Definition}

Recall that, by def\/inition, $s_{\mathbf{X}}(u) = \omega(\mathbf{X}(u))$. Also, since $h_u^t = \mathrm{exp}(tX_{h,u})$, we have $X_{h,u} = (d/dt)_{\vert t=0} h_u^t$. By def\/inition, the integral curve $\Phi_{\mathbf{X},t}(u)$ satisf\/ies $\Phi'_{\mathbf{X}}(0) = \mathbf{X}(u)$. Hence, $s_{\mathbf{X}}(u) = \omega((d/dt)_{\vert t=0}(u.h_u^t)),$ and so by property (\ref{fund-vfs}) of the Cartan connection $\omega$, we have
\[
X_{h,u} = s_{\mathbf{X}}(u).
\]
In particular, this implies the following equivariance properties of $X_{h,u}$ and $h_u^t$ with respect to a change of the base point $u \in \Gbdle_x$, so it makes sense to speak of the holonomy $h^t$ of $X$ at $x$ as a conjugacy class of one-parameter groups in $P$:
\begin{gather*}
X_{h,u.p} = \mathrm{Ad}(p^{-1})(X_{h,u}), \qquad \forall \,\, p \in P;\\
h_{u.p}^t = p^{-1}h_u^tp, \qquad \forall \,\, p \in P.
\end{gather*}

The f\/irst part of relating essentiality of $X$ near $x$ to its holonomy, is the following:

\begin{Proposition} \label{necessary holonomy} Let $X \in \mathrm{inf}(\mathcal{M})$ have a singularity $x \in M$, and let $u \in \Gbdle_x$ be as above. If $X \in \mathrm{inf}(\sigma)$ for any locally defined Weyl structure $\sigma$, then $X$ has holonomy~$h_u^t$ conjugate under~$P$ to a one-parameter subgroup of~$G_0$ $($equivalently, $s_{\mathbf{X}}(u.p) = \mathrm{Ad}(p^{-1})(s_{\mathbf{X}}(u)) \in \g_0$ for some $p \in P)$.
\end{Proposition}

\begin{proof} Assume that $\sigma$ is any locally def\/ined Weyl structure in a neighborhood of the point $x$ with $X \in \mathrm{inf}(\sigma)$. From Lemma \ref{inessential automorphisms TFAE}, it follows that $\mathbf{X}_0 \in \mathrm{inf}(\Gbdle_0,\sigma^*\omega_{\leq})$, i.e.\ we have $\Lie_{\mathbf{X}_0} \sigma^*\omega_{\leq} = 0$. Then note that the identity (\ref{infinitesimal deformation identity}) from Lemma \ref{infinitesimal deformation} simplif\/ies, for $u \in \Gbdle_x$ and $u_0 := \pi_+(u) \in (\Gbdle_0)_x$, to give us the following two identities:
\begin{gather}
(\nabla^{\Abdle} s_{\mathbf{X}})(u) = 0; \label{Have0} \\
(\nabla^{\sigma} s_{\mathbf{X}_0})(u_0) = 0. \label{Need0}
\end{gather}
 To prove the claim in the proposition, we compute the identity (\ref{Have0}) in terms of $\sigma$, to show that~(\ref{Need0}) implies $s_{\mathbf{X}}(\sigma(u_0)) \in \g_0$. Since $\sigma(u_0)$, $u \in \Gbdle_x$, therefore $s_{\mathbf{X}}(\sigma(u_0)) = \mathrm{Ad}(p^{-1})(s_{\mathbf{X}}(u)) \in \g_0$ for some $p \in P$ as claimed.

 For the computation, note that in general, a Weyl structure $\sigma$ allows us to identify a section $s \in \Gamma(\Abdle)$ with $s \circ \sigma \in C^{\infty}(\Gbdle_0,\g)^{G_0} \isom \oplus_{i=-k}^k C^{\infty}(\Gbdle_0,\g_i)^{G_0}$. We will write $[s]^{\sigma} = (s^{\sigma}_{-k}, \ldots , s^{\sigma}_k) \in \oplus_{i=-k}^k C^{\infty}(\Gbdle_0,\g_i)^{G_0}$.

  Now consider any $Y \in T_xM$, and let $\mathbf{Y}_0$ be a local right-invariant vector f\/ield on $\Gbdle_0$ around $u_0 \in (\Gbdle_0)_x$, projecting onto $Y$ at $x$, and let $\mathbf{Y}$ be a local right-invariant vector f\/ield on $\Gbdle$ which extends the vector f\/ield $\sigma_* \mathbf{Y}_0$. Then by the chain rule, we have $\mathbf{Y}(s)(\sigma(u'_0)) = \mathbf{Y}_0([s]^{\sigma})(u'_0)$ for~$u'_0$ near~$u_0$ in~$\Gbdle_0$. Now compute from the def\/inition, for $s = s_{\mathbf{X}}$ as above:
\begin{gather*}
\big(\nabla^{\Abdle}_Y s\big)(\sigma(u_0))  = \mathbf{Y}(s)(\sigma(u_0)) + \{\omega \circ \mathbf{Y},s\}(\sigma(u_0))
  = \mathbf{Y}_0([s]^{\sigma})(u_0) + \{\sigma^*\omega \circ \mathbf{Y}_0,[s]^{\sigma}\}(u_0).
\end{gather*}

 Now we translate the last line into vector notation, where the top, middle and bottom components correspond, respectively, to the projection onto $\p_+$, $\g_0$ and $\g_-$ (denoted, as usual, with a~subscript). Note that from $\Pi(s) = X$, and since $X(x) = 0$, we have $s^{\sigma}_-(\sigma(u_0)) = 0$, and so we get the following reformulation of the left-hand side of~(\ref{Have0}):
\begin{gather}
\left(\begin{array}{c}
    \mathbf{Y}_0(s^{\sigma}_+)(u_0) \\
    \mathbf{Y}_0(s^{\sigma}_0)(u_0) \\
    \mathbf{Y}_0(s^{\sigma}_-)(u_0)\end{array}\right) +
\left(\begin{array}{c}
    \{\sigma^*\omega(\mathbf{Y}_0),[s]^{\sigma}\}_+(u_0) \\
    \{\omega_0(\mathbf{Y}_0),s_0^{\sigma}\}(u_0) + \{\omega_{-}(\mathbf{Y}_0),s^{\sigma}_+\}_0(u_0) \\
    \{\omega_{-}(\mathbf{Y}_0),s^{\sigma}_0\}(u_0) +
    \{\omega_{-}(\mathbf{Y}_0),s^{\sigma}_+\}_-(u_0)
     \end{array}\right). \label{Have0 wrt sigma}
\end{gather}

On the other hand, let us compute the identity (\ref{Need0}). The section $s_{\mathbf{X}_0} \in C^{\infty}(\Gbdle_0,\p^*)^{G_0}$ is given by
\[
s_{\mathbf{X}_0}(u'_0) = \sigma^* \omega_{\leq}(\mathbf{X}_0(u'_0)) = \omega_{\leq}(\sigma(u'_0))(\sigma_*(\mathbf{X}_0(u'_0))).
\]
Using the facts that $\sigma$ is a section of $\pi_+:\Gbdle \rightarrow \Gbdle_0$, and that $\mathbf{X}$ projects onto $\mathbf{X}_0$ via $\pi_+$, it follows that $\mathbf{X}(\sigma(u'_0)) - \sigma_*(\mathbf{X}_0(u'_0))$ lies in the kernel of $T_{\sigma(u'_0)} \pi_+$. In particular, this means we have:
\[
\omega_{\leq}(\sigma(u'_0))(\sigma_*(\mathbf{X}_0(u'_0))) = \omega_{\leq}(\mathbf{X}(\sigma(u'_0))),
\]
or equivalently, $s_{\mathbf{X}_0} = s^{\sigma}_- + s^{\sigma}_0$. Using this, a similar calculation to the one above gives:
\begin{gather}
(\nabla^{\sigma}_Y s_{\mathbf{X}_0})(u_0) =
    \left(\begin{array}{c}
    \mathbf{Y}_0(s^{\sigma}_0)(u_0) \\
    \mathbf{Y}_0(s^{\sigma}_-)(u_0)\end{array}\right) +
    \left(\begin{array}{c}
    \{\omega_0(\mathbf{Y}_0),s_0^{\sigma}\}(u_0) \\
    \{\omega_{-}(\mathbf{Y}_0),s^{\sigma}_0\}(u_0) \end{array}\right).
\label{Need0 wrt sigma}
\end{gather}

 Comparing the $\g_0$-components of (\ref{Have0 wrt sigma}) and (\ref{Need0 wrt sigma}), we see that if both terms vanish, we must have
\begin{gather*}
\{\omega_{-}(\mathbf{Y}_0),s^{\sigma}_+\}_0(u_0) := \mathrm{pr}_{\g_0}([\omega_{-}(\mathbf{Y}_0)(u_0),s^{\sigma}_+(u_0)]) = 0.
\end{gather*}
But we have $\g_- = \{\omega_{-}(\mathbf{Y}_0)(u_0) \, \vert \, Y \in T_xM\}$, and from this it follows that $s^{\sigma}_+(u_0) = 0$, by using the properties of $\vert k \vert$-graded semi-simple Lie algebras: We have the grading element $E \in \g_0$, which satisf\/ies $[E,Y_j] = jY_j$ for all $Y_j \in \g_j$. Now using the $Ad$-invariance of the Killing form $B$, we have, for any $Y \in \g_-$, and $0 < j \leq k$:
\[
B([Y,s^{\sigma}_j(u_0)],E) = B(Y,[s^{\sigma}_j(u_0),E]) = -jB(Y,s^{\sigma}_j(u_0)),
\]
and since $B$ induces an isomorphism $\g_j \cong (\g_{-j})^*$, this vanishes for all $Y \in \g_-$ only if $s^{\sigma}_j(u_0) = 0$ for all $j > 0$, i.e.\ only if $s(\sigma(u_0)) \in \g_0$, which is the claim of the proposition.
\end{proof}

Next we prove the converse to Proposition~\ref{necessary holonomy}:

\begin{Proposition} \label{sufficient holonomy} Let $X \in \mathrm{inf}(\mathcal{M})$ have a singularity $x \in M$, and let $u \in \Gbdle_x$, $u_0 = \pi_+(u) \in (\Gbdle_0)_x$. If the holonomy $h_u^t$ is conjugate under $P$ to a one-parameter subgroup of $G_0$ $($equivalently, $s_{\mathbf{X}}(u.p) = \mathrm{Ad}(p^{-1})(s_{\tilde{X}}(u)) \in \g_0$ for some $p \in P)$, then $X$ is an infinitesimal automorphism of some local Weyl structure~$\sigma$ around~$x$.
\end{Proposition}

\begin{proof} We will need the following ``exponential coordinates'' on $\Gbdle$ and $M$ around $u$ and $x$ respectively, which are induced by the Cartan connection $\omega$: For any $Y \in \g$, denote by $\hat{Y}$ the vector f\/ield on $\Gbdle$ which is determined by the identity, $\omega(\hat{Y}(u')) = Y$ for all $u' \in \Gbdle$. Def\/ining
\[
\mathcal{W}_u := \{ Y \in \g \, \vert \, \varphi_{\hat{Y},t}(u) \, \mathrm{is} \, \mathrm{def\/ined} \, \mathrm{for} \, 0 \leq t \leq 1 \},
\]
then there exist an open neighborhood $\mathcal{V}_u$ of $0 \in \g$ and an open neighborhood $V_u$ of $u \in \Gbdle$ such that the exponential map $\mathrm{exp}_u^{\omega}$ def\/ined on $\mathcal{W}_u$ is a dif\/feomorphism of $\mathcal{V}_u$ onto $V_u$, where by def\/inition:
\[
\mathrm{exp}_u^{\omega}: \  Y \mapsto \mathrm{exp}^{\omega}(u,Y) := \varphi_{\hat{Y},1}(u).
\]
Restricting $\mathcal{V}_u$ if necessary to a smaller neighborhood of zero, we get a dif\/feomorphism
\[
\overline{\mathrm{exp}}_u^{\omega} := \pi \circ \mathrm{exp}_u^{\omega}: \ \mathcal{V}_u^- \stackrel{\approx}{\rightarrow} U_x,
\]
where $U_x$ is a neighborhood of $x$ in $M$ and $\mathcal{V}_u^- := \mathcal{V}_u \cap \g_-$. Furthermore, for $V_u^- := \mathrm{exp}^{\omega}_u(\mathcal{V}_u^-)$, the restriction of the projection $\pi$ gives a dif\/feomorphism of $V_u^-$ onto $U_x$.

These exponential coordinates can obviously be used to def\/ine a local Weyl structure over~$U_x$, since they give a local section of $\pi: \Gbdle \rightarrow M$ on $U_x$. This gives us the local trivialisations $\pi^{-1}(U_x) \isom V_u^- \times P$ and $\pi_0^{-1}(U_x) \isom \pi_+(V_u^-) \times G_0$. Then we simply def\/ine $\sigma: \pi_0^{-1}(U_x) \rightarrow \pi^{-1}(U_x)$ by
\begin{gather}
\sigma: \ \pi_+(u').g_0 \mapsto u'.g_0; \label{exponential Weyl structure}
\end{gather}
for any $u' \in V_u^-$, $g_0 \in G_0$. This is by def\/inition a $G_0$-equivariant local section of $\pi_+: \Gbdle \rightarrow \Gbdle_0$, that is a Weyl structure. We will now show, assuming $h_u^t \subset G_0$, that the local f\/lows of $\mathbf{X} \in \VF(\Gbdle)$ commute with $\sigma$, i.e.\ we have $\Phi_{\mathbf{X},t} \circ \sigma = \sigma \circ \Phi_{\mathbf{X}_0,t}$ on $\pi_0^{-1}(U_x)$, for $t$ suf\/f\/iciently small so that both sides exist. By Def\/inition \ref{parabolic essential definition}, this implies $\mathbf{X} \in \mathrm{inf}(\sigma)$.

To show the commutativity of the local f\/lows of $\mathbf{X}$ with the Weyl structure $\sigma$ given by~(\ref{exponential Weyl structure}), we need the following general equivariance relation for an inf\/initesimal automorphism in exponential coordinates (cf.\ the proof of Proposition~4.2 of~\cite{Frances3}):
\begin{gather}
\Phi_{\mathbf{X},t}(\mathrm{exp}^{\omega}(u,Y)) = \mathrm{exp}^{\omega}(u,\mathrm{Ad}(h_u^t)(Y)).h_u^t. \label{holonomy inf aut equiv}
\end{gather}
The identity (\ref{holonomy inf aut equiv}) is based on the observation that we have $[\mathbf{X},\hat{Y}] = 0$ for any inf\/initesimal automorphism, and any $Y \in \g$ (this follows immediately from the def\/ining equation, $\Lie_{\mathbf{X}}\omega = 0$, of an inf\/initesimal automorphism). Hence, the f\/lows commute, $\Phi_{\mathbf{X},t} \circ \Phi_{\hat{Y},s} = \Phi_{\hat{Y},s} \circ \Phi_{\mathbf{X},t}$ whenever both sides are def\/ined, which together with equivariance of $\omega$ may be used to show that both sides of (\ref{holonomy inf aut equiv}) are given as the endpoint of the same integral curve through $\Phi_{\mathbf{X},t}(u) = u.h_u^t$.

Now consider an arbitrary point $\pi_+(u').g_0 \in \pi_0^{-1}(U_x)$, where $u' = \mathrm{exp}^{\omega}(u,Y) \in V_u^-$ for $Y \in \mathcal{V}_u^-$. Then $\sigma(\pi_+(u').g_0) := u'.g_0$, and we have, by $P$-equivariance of $\Phi_{\mathbf{X},t}$ and (\ref{holonomy inf aut equiv}):
\[
(\Phi_{\mathbf{X},t} \circ \sigma)(\pi_+(u').g_0) = \Phi_{\mathbf{X},t}(u'.g_0) = \mathrm{exp}^{\omega}(u,\mathrm{Ad}(h_u^t)(Y)).h_u^tg_0.
\]
But since $h_u^t \in G_0$, we have $\mathrm{Ad}(h_u^t)(Y) \in \g_-$ and for $t$ suf\/f\/iciently small we may also assume $\mathrm{Ad}(h_u^t)(Y) \in \mathcal{V}_u$ by continuity, so $\mathrm{exp}^{\omega}(u,\mathrm{Ad}(h_u^t)(Y)) \in V_u^-$. We also have $h_u^tg_0 \in G_0$, so $G_0$-equivariance of $\pi_+$ gives $\pi_+(\Phi_{\mathbf{X},t}(u'.g_0)) = \pi_+(\mathrm{exp}^{\omega}(u,\mathrm{Ad}(h_u^t)(Y))).h_u^tg_0,$ and hence combining the above gives:
\[
(\Phi_{\mathbf{X},t} \circ \sigma)(\pi_+(u').g_0) = \Phi_{\mathbf{X},t}(u'.g_0) = (\sigma \circ \pi_+ \circ \Phi_{\mathbf{X},t})(u',g_0).
\]
Finally, since $\Phi_{\mathbf{X}_0,t}$ is def\/ined on $\Gbdle_0$ via the relation $(\Phi_{\mathbf{X}_0,t} \circ \pi_+) = (\pi_+ \circ \Phi_{\mathbf{X},t})$, this shows that $(\Phi_{\mathbf{X},t} \circ \sigma) = (\sigma \circ \Phi_{\mathbf{X}_0,t})$ on~$\pi_0^{-1}(U_x)$.
\end{proof}

\subsection{Proof of Theorem~\ref{local essential dictionary}} \label{section4.3}

There are two parts to establishing Theorem~\ref{local essential dictionary}, and they are extensions of the arguments in the proofs of Propositions~\ref{necessary holonomy} and~\ref{sufficient holonomy}.

First, we establish the following claim, which is a strengthening of Proposition~\ref{necessary holonomy}: If $X$ is inessential in some neighborhood of $x$, then its holonomy $h_u^t$ is conjugate under $P$ to a one-parameter subgroup of $\mathrm{Ker}(\lambda) \subset G_0$ for a choice of scale representation $\lambda: G_0 \rightarrow \R_+$ (equivalently, $s_{\mathbf{X}}(u.p) = \mathrm{Ad}(p^{-1})(s_{\mathbf{X}}(u)) \in \mathrm{Ker}(\lambda') \subset \g_0$ for some $p \in P$).

The proof of this claim is analogous to the proof of Proposition \ref{necessary holonomy}. If $\sigma$ is exact, then we have the holonomy reduction $\overline{\Gbdle}_0 \subset \Gbdle_0$ to structure group $\mathrm{Ker}(\lambda) \subset G_0$, and denoting the resulting reduction by $\overline{\sigma}: \overline{\Gbdle}_0 \rightarrow \Gbdle$, the condition that $X$ is an \emph{exact} inf\/initesimal automorphism of $\sigma$ is (using Lemma \ref{inessential automorphisms TFAE}), in addition to the requirements computed in the proof of Proposition \ref{necessary holonomy}, that $\mathbf{X}_0(\overline{u}) \in T_{\overline{u}}\overline{\Gbdle}_0 \subset T_{\overline{u}}\Gbdle_0$ for all $\overline{u} \in \overline{\Gbdle}_0$ and that the resulting vector f\/ield $\overline{\mathbf{X}}_0$ on $\overline{\Gbdle}_0$ is an inf\/initesimal automorphism of the Cartan connection $\overline{\sigma}^*\omega_{\leq}$. This gives
\[
\big(\nabla^{\overline{\sigma}} s_{\overline{\mathbf{X}}_0}\big)(\overline{u}_0) = 0,
\]
where $\nabla^{\overline{\sigma}}$ denotes the tractor connection on the adjoint tractor bundle associated to $\overline{\Gbdle}_0$ by the adjoint representation on $\g_- \dsum \mathrm{Ker}(\lambda')$. Then by the same considerations, if we write $s^{\sigma}_0(\overline{u}_0) = s^{\overline{\sigma}}_0(\overline{u}_0) + z(\overline{u}_0)E_{\lambda}$, for $\overline{u}_0 \in (\overline{\Gbdle}_0)_x$, then this condition implies that
\[
0 = [\omega_-(\mathbf{Y}_0(\overline{u}_0)),z(\overline{u}_0)E_{\lambda}] = \sum_{j=1}^k jz(\overline{u}_0)\omega_{-j}(\mathbf{Y}_0(u_0)),
\]
for all $Y \!\in\! T_xM$, which can only happen if $z(\overline{u}_0) = 0$. Hence we must have $s(\sigma(\overline{u}_0)) \!\in\! \mathrm{Ker}(\lambda') \!\subset\! \g_0$.

Second, we make the following, additional claim under the setting of Proposition~\ref{sufficient holonomy}: If the holonomy $h_u^t$ of $X$ is conjugate under $P$ to a one-parameter subgroup of $\mathrm{Ker}(\lambda) \subset G_0$ for a choice of scale representation $\lambda: G_0 \rightarrow \R_+$ (equivalently, $s_{\mathbf{X}}(u.p) = \mathrm{Ad}(p^{-1})(s_{\mathbf{X}}(u)) \in \mathrm{Ker}(\lambda') \subset \g_0$ for some $p \in P$), then~$X$ is inessential on some neighborhood of~$x$.

This claim is also proved in the same way as the claim of Proposition \ref{sufficient holonomy}. The considerations of that proof also show that the local, equivariant section $\sigma$ of $\pi_+: \Gbdle \rightarrow \Gbdle_0$ over $U_x$ which was constructed using the exponential map, can be restricted to a map $\overline{\sigma}: (\overline{\Gbdle}_0)_{U_x} \rightarrow \Gbdle_{U_x}$ where the sub-bundle $(\overline{\Gbdle}_0)_{U_x} \subset (\Gbdle_0)_{U_x}$ is just given by $\pi_+(V_u^-) \times \mathrm{Ker}(\lambda)$ in the local trivialization, giving a locally exact Weyl structure. Finally, since $h_u^t$ is conjugate to a one-parameter subgroup of $\mathrm{Ker}(\lambda)$, it can be arranged that the restriction of $\mathbf{X}_0$ to $(\overline{\Gbdle}_0)_{U_x}$ is always tangent to this sub-bundle, so $X$ is locally inessential by Lemma \ref{inessential automorphisms TFAE}: Namely, changing $u \in \Gbdle_x$ if necessary, we may assume that $h_u^t \subset \mathrm{Ker}(\lambda)$ and (equivalently) $\omega(\mathbf{X}(u)) \in \mathrm{Ker}(\lambda')$. For any $u' \in V_u^- := \mathrm{exp}_u(\mathcal{V}_u^-) \subset \Gbdle_{U_x}$, the identity (\ref{holonomy inf aut equiv}) can be dif\/ferentiated to show that
\[
\mathbf{X}(u') = \frac{d}{dt}_{\vert t=0}\big(\mathrm{exp}^{\omega}(u,\mathrm{Ad}(h_u^t)(Y)).h_u^t\big),
\]
for some $Y \in \mathcal{V}_u^-$. Since $h_u^t \subset \mathrm{Ker}(\lambda)$, this shows that $\mathbf{X}(u')$ is the sum of a vector tangent to~$V_u^-$ with a vertical vector corresponding to an element of $\mathrm{Ker}(\lambda') \subset \g_0$. Thus, $\mathbf{X}_0(\pi_+(u')) \in T\overline{\Gbdle}_0$ and hence this holds for restriction of $\mathbf{X}_0$ to $(\overline{\Gbdle}_0)_{\vert U_x}$, since any point in the latter space is given by $\pi_+(u').k$ for some $u' \in V_u^-$ and some $k \in \mathrm{Ker}(\lambda)$.

\subsection*{Acknowledgements} I am grateful to Charles Frances for discussions about the methods used in~\cite{Frances} and~\cite{Frances3}, and to Felipe Leitner for useful comments on an earlier version of the text. The anonymous referees made many helpful criticisms and suggestions; in particular, I am grateful for the suggestion to reformulate an earlier (equivalent) version of Def\/inition~\ref{parabolic essential definition} in terms of the action of an automorphism on the set of Weyl structures.

\pdfbookmark[1]{References}{ref}
\LastPageEnding


\begin{thebibliography}{99}

\footnotesize\itemsep=0pt


\bibitem{Aleks}
 Alekseevski D.V.,
Groups of conformal transformations of Riemannian spaces,
\href{http://dx.doi.org/10.1070/SM1972v018n02ABEH001770}{{\it Sb. Math.}} \textbf{18} (1972), 285--301.

\bibitem{Biquard}
 Biquard O.,
M\'etriques d'Einstein asymptotiquement sym\'etriques,
{\it Ast\'erisque} (2000), no.~265 (English transl.: Asymptotically symmetric Einstein metrics, {\it SMF/AMS Texts and Monographs}, Vol.~13,
 American Mathematical Society, Providence, 2006).

\bibitem{Cap inf aut}
 \v{C}ap A.,
 Inf\/initesimal automorphisms and deformations of parabolic geometries,
\href{http://dx.doi.org/10.4171/JEMS/116}{{\it J.~Eur. Math. Soc.}} \textbf{10} (2008), 415--437,
\href{http://arxiv.org/abs/math.DG/0508535}{math.DG/0508535}.

\bibitem{CS03}
\v{C}ap A., Slov\'ak J.,
Weyl structures for parabolic geometries,
{\it Math. Scand.} \textbf{93}  (2003), 53--90,
\href{http://arxiv.org/abs/math.DG/0001166}{math.DG/0001166}.

\bibitem{CSbook}
 \v{C}ap A., Slov\'ak J.,
 Parabolic geometries. I~Background and general theory,
{\it Mathematical Surveys and Monographs}, Vol.~154, American Mathematical Society, Providence, RI, 2009.

\bibitem{DuistKolk}
 Duistermaat J.J., Kolk J.A.C.,
 Lie groups, {\it Universitext}, Springer-Verlag, Berlin, 2000.

\bibitem{Ferrand}
 Ferrand J.,
 The action of conformal transformations on a Riemannian manifold,
\href{http://dx.doi.org/10.1007/BF01446294}{{\it Math. Ann.}} \textbf{304} (1996), 277--291.

\bibitem{Frances}
 Frances C.,
 Sur le groupe d'automorphismes des g\'{e}om\'{e}tries paraboliques de rang un,
\href{http://dx.doi.org/10.1016/j.ansens.2007.07.003}{{\it Ann. Sci. \'Ecole Norm. Sup.~(4)}} \textbf{40}  (2007), 741--764
 (English version:  A Ferrand--Obata theorem for rank one parabolic geometries,
\href{http://arxiv.org/abs/math.DG/0608537}{math.DG/0608537}).

\bibitem{FrEss}
 Frances C.,
 Essential conformal structures in Riemannian and Lorentzian geometry,
 in Recent Developments in Pseudo-Riemannian Geometry, {\it ESI Lectures in Mathematics and Physics}, {\it ESI Lect. Math. Phys.},
 Eur. Math. Soc., Z\"urich, 2008, 231--260.

\bibitem{Frances3}
 Frances C.,
 Local dynamics of conformal vector f\/ields,
\href{http://arxiv.org/abs/0909.0044}{arXiv:0909.0044}.

\bibitem{Gover}
Gover A.R.,
Laplacian operators and $Q$-curvature of conformally Einstein manifolds,
\href{http://dx.doi.org/10.1007/s00208-006-0004-z}{{\it Math. Ann.}} \textbf{336} (2006), 311--334,
\href{http://arxiv.org/abs/math.DG/0506037}{math.DG/0506037}.

\bibitem{IvVass}
 Ivanov S., Vassilev D.,
 Conformal quaternionic contact curvature and the local sphere theorem,
\href{http://dx.doi.org/10.1016/j.matpur.2009.11.002}{{\it J.~Math. Pures Appl.~(9)}} \textbf{93} (2010), 277--307,
\href{http://arxiv.org/abs/0707.1289}{arXiv:0707.1289}.

\bibitem{Schoen}
Schoen R.,
On the conformal and CR automorphism groups,
\href{http://dx.doi.org/10.1007/BF01895676}{{\it Geom. Funct. Anal.}} \textbf{5} (1995),   464--481.

\bibitem{Webster}
Webster S.M.,
On the transformation group of a real hypersurface,
\href{http://dx.doi.org/10.2307/1997877}{{\it Trans. Amer. Math. Soc.}} \textbf{231} (1977), 179--190.

\end{thebibliography}
\end{document}